\documentclass[12pt,a4paper]{article}
\usepackage{amsmath,amssymb,amsfonts}
\def\Bbb{\mathbb}

\title{\bf On the equality of Dedekind sums}

\author{Kurt Girstmair\\
Institut f\"ur Mathematik, Universit\"at Innsbruck   \\
Technikerstr. 13/7, A-6020 Innsbruck, Austria \\
Kurt.Girstmair@uibk.ac.at}
\date{}

\makeatletter
\let\@@maketitle=\maketitle
\def\maketitle{\def\thispagestyle##1{\relax}\@@maketitle}
\makeatother
\newtheorem{theorem}{Theorem}
\newtheorem{prop}{Proposition}
\newtheorem{lemma}{Lemma}
\newtheorem{corollary}{Corollary}

\def\BE{\begin{equation}}
\def\EE{\end{equation}}
\def\BD{\begin{displaymath}}
\def\ED{\end{displaymath}}
\def\BA{\begin{array}}
\def\EA{\end{array}}
\def\BEA{\begin{eqnarray*}}
\def\EEA{\end{eqnarray*}}
\def\BI{\bibitem}

\def\Z{\Bbb Z}

\def\R{\Bbb R}

\def\phi{\varphi}

\def\MB{\mbox}
\def\LD{\ldots}

\def\DIV{\,|\,}

\def\MN{\medskip\noindent}
\def\STOP{\hfill$\Box$}

\def\DED{Dedekind }
\def\CD{\{c,d\}}

\begin{document}
\maketitle

\begin{abstract}

\noindent
We show that deciding the equality of two Dedekind sums $S(c,b)$, $S(d,b)$ is equivalent to deciding whether a Dedekind sum  defined by $b, c, d$ takes a certain value. By means of this result we construct
infinite sequences of pairwise equal Dedekind sums. Moreover, we prove a result that says how many Dedekind sums $S(d,b)$, $1\le d\le b-1$, may be equal to a given $S(c,b)$ if $b$ is a square-free number.

\end{abstract}

\section*{1. Introduction}

Let $b$ be a natural number and $c$ an integer such that $(c,b)=1$. The classical \DED sum $s(c,b)$ is defined by
\BE
\label{1.0}
   s(c,b)=\sum_{k=1}^{b} ((k/b))((ck/b)).
\EE
Here
\BD
  ((x))=\begin{cases}
                 x-\lfloor x\rfloor-1/2 & \MB{ if } x\in\R\smallsetminus \Z; \\
                 0 & \MB{ if } x\in \Z
        \end{cases}
\ED
(see \cite[p. 1]{RaGr}).

 \DED sums first appeared in the theory of modular forms; see \cite{Ap2}. But these sums have also interesting applications
in a number of other fields, so in connection with class numbers, lattice point problems, topology, and algebraic geometry (see, for instance,  \cite{{At}, {BoLi}, {Me}, {RaGr}, {Ur}, {Ur2}}).

It is often more convenient to work with
\BD
 S(c,b)=12s(c,b).
\ED
We call $S(c,b)$ a {\em normalized} \DED sum. The number $b$ is called the {\em modulus} and $c$ the {\em argument} of this \DED sum.

The present paper is devoted to the study of equal values of \DED sums belonging to the same modulus $b$ but to different arguments $c,d$.
Two cases of equality can be considered as trivial.
First,
\BD
  S(c,b)=S(c',b)
\ED
 if $c\equiv c' \mod b$ (this justifies calling $b$ the modulus of $S(c,b)$).
Second, let $c^*$ denote an inverse of $c \mod b$, i.e., an the integer such that
$cc^*\equiv 1 \mod b$. Of course, $c^*$ is determined only up to the addition of multiples of $b$. It is an easy consequence of the definition (\ref{1.0}) that
\BD
  S(c,b)=S(c^*,b);
\ED
see \cite[p. 26]{RaGr}.

Cases of nontrivial equality, i. e., $S(c,b)=S(d,b)$ and $d\not\equiv c,c^*\mod b$, do not occur frequently.  For example, for $b=3^3\cdot 11=297$, $S(c,b)$ takes 41 distinct positive values,
but only five of them give rise to nontrivial equality, namely, $3076/b$, $1712/b$, $1460/b$, $1456/b$, and $1136/b$.
If $b=p_1^{m_1}\cdots p_r^{m_r}$, the $p_i$ being distinct primes, then the number of cases of
nontrivial equality seems to increase with $r$ and the exponents $m_i$.

Suppose that the modulus $b$ and two arguments $c,d$ are given, $d\not\equiv c,c^*\mod b$. So far no simple condition is known that is equivalent to $S(c,b)=S(d,b)$. One has only necessary
conditions for this equality like
\BE
\label{1.2}
     (c-d)(cd-1)\equiv 0\mod b;
\EE
see \cite{Ja}. Indeed, the condition (\ref{1.2}) is equivalent to $S(c,b)\equiv S(d,b)\mod\Z$; see \cite{Gi1}. This criterion was extended to a condition for $S(c,b)\equiv S(d,b)\mod 8\Z$; see \cite{Ts}.
However, the latter condition is no more quite simple, because it may involve up to two Jacobi symbols that must be evaluated.

In the present paper we prove the following theorem, which is based on an adaption of a result of U. Dieter; see Section 2.

\begin{theorem} 
\label{t1}
Let $c,d\in\Z$, $(c,b)=(d,b)=1$, $c\not\equiv d\mod b$. Suppose that $b,c,d$ satisfy condition {\rm(\ref{1.2})}. Let $t>0$ be such that $t\equiv c-d\mod b$. Then $S(c,b)=S(d,b)$ if, and only if,
\BD
\label{1.4}
      S(1+ct,bt)=\frac{2+t^2}{bt}-3.
\ED
\end{theorem} 

\MN
If we want to decide whether two \DED sums with the same modulus are equal, we may, instead, decide whether a certain \DED sum takes a certain value, as the theorem says.
We think that the latter task is not simpler in most cases. But, due to this theorem, it seems to be plausible that no simple necessary and sufficient condition for nontrivial equality exists.

In this paper we present two applications of Theorem \ref{t1}.
The first application is the construction of infinite sequences of pairwise equal \DED sums.
To this end let $c,d\in \Z$, $c\not\equiv d\mod b$, with $(c,b)=(d,b)=1$ (in particular, $b>1$).
We say that $\CD$ is a {\em suitable set} for (the modulus) $b$ if, and only if,
$d\not\equiv c^* \mod b$ and $S(c,b)=S(d,b)\ne 0$.

\begin{theorem} 
\label{t2}
In the above setting, let $\CD$ be a suitable set for $b$. Let $t>0$ be such that $t\equiv d-c^* \mod b$.  Put $b_1=bt$,
$c_1=1+ct$, $d_1=1+dt$. Then $b_1>b$ and $\{c_1,d_1\}$ is a suitable set for $b_1$.
\end{theorem} 

By Theorem \ref{t2}, we can define an infinite sequence of strictly increasing moduli $b_i$, $i\ge 0$, together with sets $\{c_i,d_i\}$ suitable for $b_i$, provided that one suitable set $\CD$ is known.
Indeed, put $b_0=b$, $c_0=c$, $d_0=d$ and $t_0=t$. Define, recursively, $b_{i+1}=b_it_i$, $c_{i+1}=1+c_it_i$, $d_{i+1}=1+d_it_i$, $i\ge 0$. Then define $t_{i+1}>0$
by $t_{i+1}\equiv d_{i+1}-c_{i+1}^* \mod b_{i+1}$ (the inverse is to be understood$\mod b_{i+1}$).
In particular, we obtain $d_i\not\equiv c_i, c_i^*\mod b_i$ and
\BD
 S(c_i,b_i)=S(d_i,b_i)\ne 0
\ED
for all $i\ge 0$.

\MN
{\em Example.} Put $b=7\cdot 11=77$, $c=16$, $d=60$. Then $d\not\equiv c,c^*\mod b$ and $S(c,b)=S(d,b)=300/77$. Obviously, the set $\CD$ is suitable for $b$. Since $d-c^*\equiv 7$ mod $b$
we may take $t=7$ and form the sequences $b_i$ and $\{c_i,d_i\}$, $i\ge 0$, in the above way, where $t_i$ is chosen in $\{1,\LD,b_i-1\}$.
We obtain
\BEA
&&b_0=77, c_0=16, d_0=60, \\
&&b_1=539, c_1=113, d_1=421,\\
&&b_2=260337, c_2=54580, d_2=203344,  \\
&&b_3=6412881321, c_3=1344469141, d_3=5008972753, \\
&&b_4=36852630635308805163, c_4=7726203273338872624, \\
&&d_4=28784849350658189860.
\EEA
One sees that these numbers grow rapidly. The first values of $t_i$ are $t_0=7$, $t_1=483$, $t_2=24633$, $t_3=5746657203$.

\MN
How can we find suitable sets $\{c,d\}$ to initiate sequences like the above? A partial answer is given by the following theorem, whose proof involves another application of Theorem \ref{t1}.

\begin{theorem} 
\label{t3}
Let $1\le k\le r$ and $p_1,\LD,p_k$ be distinct primes, each of which is congruent $\pm 1\mod 5$. Put $b_0=p_1\cdots p_k$ and let $p_{k+1},\LD,p_r$ be distinct primes, each of which is congruent $1\mod b_0$.
Let $b=p_1\cdots p_r$ and $t=p_{k+1}\cdots p_r$.
Then
\BD
  \left|\left\{c: 1\le c\le b-1, (c,b)=1, S(c,b)=\frac{t^2+2}{b}-3\right\}\right|=2^k.
\ED
\end{theorem} 

\MN
On observing that $t=1$ if $k=r$, we have the following.

\begin{corollary}
\label{c1}
Let $1\le r$ and $p_1,\LD,p_r$ be distinct primes, each of which is congruent $\pm 1\mod 5$. Let $b=p_1\cdots p_r$. Then
\BD
  \left|\left\{c: 1\le c\le b-1, (c,b)=1, S(c,b)=\frac{3}{b}-3\right\}\right|=2^r.
\ED
\end{corollary}

\MN
It is easy to see that the common value $(t^2+2)/b-3$ of the \DED sums in Theorem \ref{t3} cannot vanish; see the end of the proof of Theorem \ref{t2} in Section 3.

The proof of Theorem \ref{t3} shows how to find the numbers $c$ in question by means of the Chinese remainder theorem; see Section 3.
Suppose, for a moment, that $b$ is as in Corollary \ref{c1} with $r=3$. Then we have eight numbers $c$ such that $S(c,b)=3/b-3$. These numbers give us $24$ suitable sets for the modulus $b$.
Section 3 contains additional examples of suitable sets and a few remarks on the above sequences of suitable sets.

Given $b$ and $c$, $(c,b)=1$, we consider the number
\BE
\label{1.8}
  N(c,b)=|\{d: 1\le d\le b-1, (d,b)=1, S(d,b)=S(c,b)\}|.
\EE
Suppose that $b$ is a square-free number {\em consisting of $r$ primes}, i.e., $b=p_1\cdots p_r$, the $p_i$ being distinct. It is known that $N(c,b)\le 2^r$; see \cite[Th. 3]{Gi1}.
Theorem \ref{t3} exhibits the $2$-powers $2^k$, $1\le k\le r$, as possible values of $N(c,b)$ for this case.
At the end of Section 3 we will see that there may be values $>1$ of $N(c,b)$ different from the aforesaid $2$-powers for a
number $b$ of this kind.

\section*{2. The criterion}

The proof of Theorem \ref{t1} is based on the following proposition.

\begin{prop} 
\label{p1}
Let $c,d\in\Z$, $(c,b)=(d,b)=1$, $c\not\equiv d\mod b$.  Let $t>0$ be such that $t\equiv c-d\mod b$.
Then
\BE
\label{2.0}
S(1+d^*t,bt)-\left(\frac{t^2+2}{bt}-3\right) =S(d,b)-S(c,b).
\EE
\end{prop} 

\MN
Of course, the identity (\ref{2.0}) may also serve as a criterion for the equality of $S(c,b)$ and $S(d,b)$. It is, however, less simple than Theorem \ref{t1} since it involves the inversion of $d\mod b$.
This inversion involves more work than checking the condition (\ref{1.2}) (as required by Theorem \ref{t1}).

Proposition \ref{p1} can be found, in a rather disguised form, in a paper of U. Dieter; see \cite[Satz 4]{Di}. In particular, it is not obvious that Dieter's version contains a criterion
for the equality of \DED sums. Dieter obtained his result as an application of his three-term relation; see \cite[Satz 1]{Di}.
We prefer deriving Proposition \ref{p1} directly from this well-known relation since adapting Dieter's Satz 4 would not be simpler.

\MN
{\em Proof of Proposition \ref{p1}.}
The three-term relation, in its most convenient form for the present purpose, reads as follows; see \cite{Gi}. Let $B$ and $D$ be natural numbers, $A$ and $C$ integers with $(A,B)=(C,D)=1$.
Suppose that
\BD
   Q=AD-BC>0.
\ED
Let $j$ and $k$ be integers such that
\BD
 -Cj+Dk=1.
\ED
Define $R$ by
\BD
 R=Aj-Bk.
\ED
Then
\BE
\label{2.8}
 S(A,B)=S(C,D)+S(R,Q)+\frac{B^2+D^2+Q^2}{BDQ}-3.
\EE
Let $b$, $c$, $d$, and $t$ be as in the proposition.  Let $c'\equiv c\mod b$ be such that $c'-d=t$. We put $B=D=b$,
$A=c'$, and $C=d$. Then $Q=bt>0$. Since $-Cj+Dk=-dj+bk=1$, $j=-d^*$ for an inverse $d^*$ of $d \mod b$ and $k=(1-dd^*)/b$. We obtain  $R=-1-d^*t$. Since $-S(R,Q)=S(1+d^*t,bt)$ and
$(B^2+D^2+Q^2)/(BDQ)=(t^2+2)/(bt)$, the identity (\ref{2.8}) gives (\ref{2.0}), but with $c'$ instead of $c$. However, $S(c',b)=S(c,b)$, and the other quantities in (\ref{2.0}) depend only
of $t$ and $b$. Hence (\ref{2.0}) holds in the above form.
\STOP

\begin{lemma} 
\label{l1}
Suppose that $b$, $c$, $d$ satisfy {\rm(\ref{1.2})}. Then
\BD
\label{2.10}
  c-d\equiv d^*-c^*\mod b.
\ED
\end{lemma} 

\MN
{\em Proof.}
Indeed, suppose $b=p_1^{m_1}\cdots p_r^{m_r}$, where the $p_i$ are distinct primes. It suffices to show
\BD
  c-d\equiv d^*-c^*\mod p_i^{m_i} \MB{ for all } i\in\{1,\LD,r\}.
\ED
We fix $i$ for the time being. So we write $p=p_i$, $m=m_i$. The congruence (\ref{1.2}) implies
\BD
 (c-d)(cd-1)\equiv 0 \mod p^{m}.
\ED
If we multiply this congruence by $c^*$, we obtain
\BD
 (1-c^*d)(d-c^*)\equiv 0\mod p^m.
\ED
Accordingly, $1-c^*d\equiv 0\mod p^j$, $d-c^*\equiv 0\mod p^k$, where the nonnegative integers $j,k$ are such that $j+k\ge m$. Since $c^*\equiv d\mod p^k$, we have $c\equiv d^*\mod p^k$.
If we write $c=d^*+up^k$, $1-c^*d=vp^j$, $u,v\in\Z$, we obtain
\BD
 c(1-c^*d)=d^*(1-c^*d)+uvp^{j+k}
\ED
and $c-d\equiv d^*-c^*\mod p^{j+k}$. \STOP

\MN
{\em Proof of Theorem \ref{t1}.} Let $t>0$, $t\equiv c-d\mod b$. Because $b$, $c$, $d$ satisfy (\ref{1.2}), Lemma \ref{l1} yields $t\equiv d^*-c^*\mod b$.
We replace, in the setting of Proposition \ref{p1}, the number $c$ by $d^*$ and $d$ by $c^*$.
This gives
\BD
  S(1+ct,bt)-\left(\frac{t^2+2}{bt}-3\right) =S(c^*,b)-S(d^*,b)=S(c,b)-S(d,b),
\ED
whence the assertion follows.
\STOP

\section*{3. Suitable sets}

\noindent
{\em Proof of Theorem \ref{t2}.}
Let $t>0$ be as in Theorem \ref{t2}, i.e., $t\equiv c-d^*\mod b$. Since $S(c,b)=S(d^*,b)$, the numbers $c$ and $d^*$ satisfy (\ref{1.2}),
and Theorem \ref{t1} gives
\BE
\label{3.2}
  S(1+ct,bt)=\frac{t^2+2}{bt}-3.
\EE
By Lemma \ref{l1}, we also have $t\equiv d-c^*\mod b$. Since $S(d,b)=S(c^*,b)$, Theorem \ref{t1}
yields
\BD
  S(1+dt,bt)=\frac{t^2+2}{bt}-3.
\ED
Accordingly, $S(c_1,b_1)=S(d_1,b_1)$ for $b_1=bt$, $c_1=1+ct$, and $d_1=1+dt$.

It remains to be shown that
$d_1\not\equiv c_1, c_1^*\mod b_1$, $t>1$, and $S(c_1,b_1)\ne 0$.

First we check
$d_1\not\equiv c_1\mod b_1$. Since $d-c\not\equiv 0 \mod b$, we have $(d-c)t\not\equiv 0 \mod bt$.
However,
\BD
  d_1-c_1\equiv (d-c)t \mod bt,
\ED
whence the assertion follows.

We also have to exclude that $d_1$ is an inverse of $c_1\mod b_1$. From
\BD
   (1+ct)(1+dt)\equiv 1\mod bt
\ED
we obtain
\BD
  cdt^2+(c+d)t\equiv 0\mod bt \:\MB{ and }\: cdt+c+d\equiv 0\mod b.
\ED
If we use $t\equiv c-d^*\mod b$ in the last congruence, we obtain
$c^2d+d\equiv 0 \mod b$,
and so $c^2\equiv-1\mod b$. But this implies $S(c,b)=0$ (see \cite[Satz 1]{Ra}), which we excluded.

Now we show $(t,b)>1$.
The congruences $c-d^*\equiv t\mod b$ and $d-c^*\equiv t\mod b$ imply $cd-1\equiv td\mod b$ and $cd-1\equiv tc\mod b$. In particular, $td\equiv tc\mod b$.
Thus, if $(t,b)=1$, we have $c\equiv d\mod b$, which we excluded. In particular, $t=1$ is impossible.

If $S(c_1,b_1)=0$, then $t^2+2-3bt=0$, by (\ref{3.2}).
Therefore, $t\DIV 2$, and so $t=2$, since $t=1$ is excluded. If $t=2$, we have $6-6b=0$, which implies
$b=1$. But in this case a suitable set $\{c,d\}$ does not exist.
\STOP

\MN
{\em Remarks.} 1. Recall the definition of the sequences $b_i$, $\{c_i, d_i\}$, $i\ge 0$, of Section 1. The congruence $t_{i+1}\equiv d_{i+1}-c_{i+1}^*\mod b_{i+1}$ implies
$t_{i+1}\equiv d_{i+1}-c_{i+1}^*\mod t_i$. Since $c_{i+1}\equiv d_{i+1}\equiv 1\mod t_i$, this gives $t_{i+1}\equiv 0\mod t_i$ for all $i\ge 0$. So the numbers $t_i$ form an ascending chain
$t_0\DIV t_1\DIV t_2\cdots$ of divisors.

2. If we are looking for a suitable set in order to start a sequence of this type, we will be successful, as it seems, if we restrict our search to square-free numbers $b\ge 70$ consisting of exactly two
primes $\ge 5$. We present a small table of such numbers $b$ together with suitable sets $\CD $. In all cases the \DED sum belonging to the suitable set is positive.
\BD
\begin{array}{r|r|r|r|r|r|r|r|r}
b & 77 & 85 & 91 & 95 & 115 & 119 & 133 & 143\\ \hline \rule[5mm]{0mm}{1mm}
c &  9 &  7 &  5 & 33 &  18 &  31 &  54 &   8\\ \hline \rule[5mm]{0mm}{1mm}
d & 16 & 22 & 31 & 52 &  78 &  45 &  73 &  73
\end{array}
\ED

3. Suppose that we restrict $t_i$ to the range $1\le t_i\le b_i$ in the above sequence.
It would be interesting to understand the limiting behaviour of the sequence $S(c_i,b_i)$, which, by (\ref{3.2}), is equivalent to the behaviour of $t_i/b_i$.

4. Note that every set $\CD$ suitable for $b$ defines three additional suitable sets, namely, $\{c,d^*\}$, $\{c^*,d\}$, and $\{c^*, d^*\}$.  Each of these four sets defines an appropriate
number $t$. Even if we restrict $t$ to the range $1\le t\le b$, we obtain, as a rule, four possibilities for $\CD$ and $t$. Of course, one may use these possibilities
for the construction of $\{c_i, d_i\}$ and $t_i$ in each step $i\ge 0$. In this way one obtains, instead of an infinite sequence, an infinite cascade of pairwise equal \DED sums.

\MN
{\em Proof of Theorem \ref{t3}.} In the setting of this theorem,
let $i\in\{1,\LD, k\}$. Since $5$ is a quadratic residue$\mod p_i$, there is an integer $\alpha_i$ such that $\alpha_i^2\equiv 5\mod p_i$.
Define $c\in\{1,\LD, b\}$ by
\BD
  c\equiv(3+\alpha_i)/2 \mod p_i,\: i=1,\LD,k,
\ED
and $c\equiv 1\mod p_i$, $i=k+1,\LD,r$ (here $1/2$ stands for an inverse of $2$ mod $p_i$). Let $d\in\{1,\LD, b\}$ be a solution of the congruence (\ref{1.2}). This means that
\BD
  (c-d)(cd-1)\equiv 0 \mod p_i,
\ED
for all $i\in\{1,\LD,r\}$. If $i\in\{k+1,\LD,r\}$, this congruence is equivalent to $(1-d)(1-d)\equiv 0\mod p_i$, i.e., $d\equiv 1\equiv c\mod p_i$. If $i\in\{1,\LD,k\}$, $d$ must satisfy one of the congruences
\BD
  d\equiv (3+\alpha_i)/2\mod p_i \:\MB{ or }\: d\equiv (3-\alpha_i)/2\mod p_i,
\ED
since $(3-\alpha_i)/2$ is an inverse of $(3+\alpha_i)/2\mod p_i$. Altogether, $d$ can be defined by the congruences
\BD
   d\equiv(3+(-1)^{j_i}\alpha_i)/2\mod p_i, i=1,\LD,k, \:\MB { and }\: d\equiv 1\mod p_i, i=k+1,\LD,r,
\ED
where $j_i\in\{0,1\}$ may be arbitrary for each $i\in\{1,\LD,k\}$. Altogether, we have exactly $2^r$ distinct numbers $d\in\{1,\LD,b-1\}$ such that (\ref{1.2}) holds. Observe that $(c,b)=1$ since an inverse $c^*$ of $c\mod b$
is given by $c^*\equiv(3-\alpha_i)/2\mod p_1$, $i=1,\LD,k$, and $c^*\equiv 1\mod p_i$, $i=k+1,\LD,r$. Because (\ref{1.2}) is the same as saying $S(c,b)-S(d,b)\in\Z$, there are exactly
$2^k$ integers $d\in\{1,\LD,b-1\}$ such that $S(c,b)-S(d,b)\in \Z$.

Recall that $t=p_{k+1}\cdots p_r$.
We have to show that $S(c,b)=S(d,b)=(t^2+2)/b-3$ for all these integers $d$. Given such an integer $d$, we define $m\in\{1,\LD,b_0\}$ by $m\equiv d-1\mod p_i$, $i=1,\LD,k$. In other words,
$m\equiv (1+(-1)^{j_i}\alpha_i)/2\mod p_i $ for these
numbers $i$. Then $m$ is invertible mod $b_0$, an inverse $m^*$ being defined by $m^*\equiv((-1)^{j_i}\alpha_i-1)/2\mod p_i$ for these $i$. In particular, $m-m^*\equiv 1\mod b_0$.
Now $t\equiv 1\equiv m-m^*\mod b_0$.
Since $S(m, b_0)=S(m^*,b_0)$, Theorem \ref{t1} can be applied to $b_0$, $m$, and $m^*$. It gives
\BD
  S(1+mt,b_0t)=\frac{t^2+2}{b_0t}-3.
\ED
Here $1+mt\equiv 1+m\equiv d\mod p_i$, $i=1,\LD, k$, and $1+mt\equiv 1\equiv d\mod p_i$, $i=k+1,\LD,r$, because $p_i\DIV t$ for these $i$. In other words, $1+mt\equiv d\mod b$, and so $S(d,b)=(t^2+2)/b-3$.
In particular, $S(c,b)$ also takes this value, since this case corresponds to $j_1=\cdots=j_k=0$.
\MB{ }\STOP

\MN
{\em Remark.} Let $b$ be a square-free number consisting of $r$ primes. For an integer $c$ with $(c,d)=1$ let $N(c,b)$ be defined as in (\ref{1.8}). Let $r=3$. Hence $N(c,b)\le 8$, as we said in
Section 1. Theorem \ref{t3} says that $2$, $4$, and $8$ are
possible values $>1$ of $N(c,b)$. It is not difficult to see that, for a square-free number $b$, $N(c,b)$ either equals $1$ or is even. Therefore, the only possible value $>1$ not in this
list is $N(c,b)=6$. In the case $b=455=5\cdot 7\cdot 13$ we find $N(c,b)=6$ for $c=32$. On the other hand,
there are no integers $c$ with $N(c,b)=8$ in this case.


\end{document}